# Feasibility Analysis and Proof of Model Predictive Control

**Hanyu Zhang**
Department of Civil and Coastal Engineering,
University of Florida, Gainesville, Florida, 32608
Email: hanyu.zhang@ufl.edu
ORCID:0000-0002-5644-612X

**Lili Du* (Corresponding author)**
Associate professor
Department of Civil and Coastal Engineering,
University of Florida, Gainesville, Florida, 32608
Email: lilidu@ufl.edu
ORCID:0000-0003-1740-1209

**Abstract**

Sequential and terminal constraint feasibility of the model predictive control (MPC) play important roles in ensuring MPC control continuity. This study thus investigates these two properties theoretically using an MPC model for vehicle platooning and eco-driving strategy at signalized intersections. This study is served as a supplement material for our paper coming soon.

1. **Introduction**

Model predictive control (MPC) [1] is an advanced control approach for a dynamic process while satisfying a set of constraints. It iteratively predicts the future environment and dynamic states, solves an optimization model with finite prediction horizons and only implement the first control law at each control time step. Increasing research [2][3][4] have been on MPC because it is robust to the environment and control uncertainties. According to [5], the MPC sequential feasibility, also named as recursive feasibility is crucial to guarantee the MPC control continuity. It is an important theoretical property of MPC that we cannot overlook. Besides, existing MPC controllers [6] often used terminal constraints to ensure the MPC control stability. Nevertheless, the extra terminal constraints pose strict requirements to the prediction horizon for feasibility. Motived by this view, this study first analyzes the MPC sequential feasibility and then based on that proves the MPC terminal constraints feasibility of an existing MPC controller.

The roadmap of this study is as follows. Following the introduction in section 1, we will present the MPC control model in section 2. Then taking this model as example, we analyze and prove the MPC sequential and terminal constraints feasibility. Finally, we conclude this study in section 4.

2. **MPC Model**

This study considers the MPC controller below in Equations (1)-(14). Mainly, the MPC controller generates the platoon-centered trajectory control laws to guide the car-following movement of CAVs in the platoon at time step $k \in \mathbb{Z}_+$ by predicting the platoon future states at any time step $k + p$, for $\forall p = 1, \dots, P$ $(p \in P)$, before it enters the traffic signal communication zone. It aims to minimize traffic oscillations and energy consumption represented by Equation (1), subject to the vehicle dynamics and various constraints demonstrated in Equations (2)-(14), where the time steps in the prediction horizon is denoted as $k_P$. Hereafter, we also simplify step $k + p$ to $p$ throughout this section to avoid complex notation.

**MPC-$q_0$**

$$\textbf{Min } \Gamma(u(p)) = \sum_{p=1}^{P} \{\frac{1}{2}[z^T(p)Q_z z(p) + (z'(p))^T Q_{z'} z'(p)] + \frac{\tau^2}{2}\omega_1 \|u(p-1)\|_2^2\} \quad (1)$$

Subject to

$$x_i(k+1) = x_i(k) + \tau v_i(k) + \frac{\tau^2}{2}\big(u_i(k) - \Delta u_i(k)\big), \quad i \in I_C, k \in k_p \tag{2}$$

$$v_i(k+1) = v_i(k) + \tau\big(u_i(k) - \Delta u_i(k)\big), \quad i \in I_C, k \in k_p \tag{3}$$

$$\Delta u_i(k) = \varepsilon_i v_i(k) + \eta_i u_i(k) - \eta_i u_i(k-1), \quad i \in I_C, k \in k_p \tag{4}$$

$$x_{\widehat{m}}(k) = x_n(k - T_{\widehat{m}}) - D_{\widehat{m}}, \quad k \in k_p \tag{5}$$

$$a_{min,i} \leq u_i(k) \leq a_{max,i}, \quad i \in I_C, k \in k_p \tag{6}$$

$$v_{min} \leq v_i(k) \leq v_{max}, \quad i \in I_C, k \in k_p \tag{7}$$

$$x_{i-1}(k) - x_i(k) \geq L_i + \delta_1 \tau v_i(k) + \delta_2 \tau\big(v_i(k) - v_{i-1}(k)\big), \quad i \in I_C, k \in k_p \tag{8}$$

$$s_i(k) = L_i + \delta_1 \tau v_i(k) + \delta_2 \tau\big(v_i(k) - v_{i-1}(k)\big) + \delta, \quad i \in I_C, k \in k_p \tag{9}$$

$$\Delta x_i(k) = x_{i-1}(k) - x_i(k) - s_i(k), \quad i \in I_C, k \in k_p \tag{10}$$

$$\Delta v_i(k) = v_{i-1}(k) - v_i(k), \quad i \in I_C, k \in k_p \tag{11}$$

$$z(k) \coloneqq \big(\Delta x_1(k), \dots, \Delta x_N(k)\big)^T \in \mathbb{R}^N, \quad k \in k_p \tag{12}$$

$$z'(k) \coloneqq \big(\Delta v_1(k), \dots, \Delta v_N(k)\big)^T \in \mathbb{R}^N, \quad k \in k_p \tag{13}$$

$$\Phi_f \coloneqq (z(k+P) \in \zeta, z'(k+P) \in \zeta'). \tag{14}$$

More exactly, Equations (2)-(4) represent the robust CAV dynamics using double integrator model. $\Delta u_i(k)$ represents the CAV $i$'s control uncertainties including powertrain delay and aerodynamic drag etc. Equation (5) uses Newell's car-following model (Newell, 2002) to predict HDV driving behaviors. Equations (6)-(8) respectively illustrate the acceleration, speed limits and safe distance constraints of CAV $i$. Equation (9) describes the desired spacing policy for CAV $i$, Equations (10)-(11) curve the spacing and speed tracking errors of CAV $i$. Accordingly, the tracking errors of all the CAVs together form the platoon tracking dynamics in Equations (12)-(13). Finally, Equation (14) presents the terminal constraints. It requires the platoon spacing and speed tracking errors will be confined to small domains $\zeta$ and $\zeta'$ respectively at final time step $k+P$ of the MPC prediction horizon.

## 3. MPC Terminal Constraints Feasibility

This section analyzes and proves the sequential feasibility and terminal constraint feasibility of the MPC controller above. They together ensure the control continuity and smoothness of the MPC system.

### 3.1. MPC sequential feasibility

MPC is implemented recursively at each time step $0, 1, \dots, k-1, k$. Therefore, a fundamental theoretical question is whether the MPC can find a feasible control law at each time step $k$ (i.e., whether the constraint set of the MPC optimizer is non-empty at each time step $k$), given the platoon system starts from an initial feasible condition at $k=0$. A MPC system is called sequential (recursive) feasible (Löfberg, 2012) if the answer to this question is affirmative. The hybrid MPC system in this study has three controllers: MPC-$q_0$, MPC-$q_1$ and MPC-$q_2$. MPC-$q_0$ has constraints in Equations (2)-(14), which are also shared with MPC-$q_1$ and MPC-$q_2$ except the terminal constraint in Equation (14). Hence, this study first proves the sequential feasibility of MPC-$q_0$. Then, we further discuss the sequential feasibility of MPC-$q_1$ and MPC-$q_2$ as well as the switching feasibility of the hybrid system in section 4.2. To prove the sequential feasibility of the MPC-$q_0$, we first classify Equations (2)-(14) into the following three sets:

(i) $S_1(u(k))$: constraint set in Equations (2)-(4) and (6)-(8) for capturing the CAV dynamics, acceleration, speed and safety constraints at step $k \in \mathbb{Z}_+$.

(ii) $S_2(u(k), Z(k+P))$: the HDV movements in Equation (5) and the terminal constraint in Equation (14) at step $k \in \mathbb{Z}_+$.

(iii) $S_3(u(k))$: the control dynamics in Equations (9)-(13) at step $k \in \mathbb{Z}_+$.

It should be noticed that the third constraint set $S_3(u(k))$ involves control dynamic formulations. They are always feasible if the first constraint set $S_1(u(k))$ is feasible. For the second constraint set $S_2(u(k), Z(k+P))$, Equations (5) is an equality constraint to curve the HDV trajectory. Thus, they are always feasible in math and also stay feasible in practice if accurate time and distance displacements are estimated. Apart from it, the terminal constraint in Equation (14) is only active at the final time step $k+P$ of the MPC. Its sequential feasibility is ensured if the nominal MPC system is asymptotical stable (see proof in detail in section 5.2) according to Mayne et al., (2000). Consequently, to prove the sequential feasibility of the MPC-$q_0$, this study will mainly analyze and prove the sequential feasibility of the first constraint set $S_1(u(k))$ by **Lemma 1**.

**Lemma 1.** For $k \in \mathbb{Z}_+ := \{0, 1, 2, \ldots\}$ and $i \in I_C$, if $S_1(u_i(k))$ is feasible, then there exists $\delta_1 \geq 1$ and $\delta_2 \geq 0$ that make $S_1(u_i(k+1))$ feasible and compact. In addition, the non-empty feasible control input profile $\mathbb{S}_1(u_i(k))$ for platoon vehicle $i$ at step $k$ is given below:

$$u_i(k) \in \mathbb{S}_1(u_i(k)) = \left[\max\{a_{min,i}, \underline{a_{i,v}}\}, \min\{a_{max,i}, \overline{a_{i,v}}, \overline{a_{i,d}}\}\right], \tag{15}$$

where

$$\underline{a_{i,v}} = \frac{v_{min} - (1 - \tau\varepsilon_i)v_i(k) - \tau\eta_i u_i(k-1)}{\tau(1 - \eta_i)} \leq 0$$

$$\overline{a_{i,v}} = \frac{v_{max} - (1 - \tau\varepsilon_i)v_i(k) - \tau\eta_i u_i(k-1)}{\tau(1 - \eta_i)} \geq 0$$

$$\overline{a_{i,d}} = \varepsilon_i v_i(k) + \eta_i u_i(k) - \eta_i u_i(k-1) + \frac{g_i(k) + \tau(v_{i-1}(k+1) - v_i(k))}{\tau^2\left(\delta_1 + \delta_2 + \frac{1}{2}\right)}$$

$$+ \frac{\left(\delta_2 - \frac{1}{2}\right)u_{i-1}(k)}{\delta_1 + \delta_2 + \frac{1}{2}}$$

$$g_i(k) = x_{i-1}(k) - x_i(k) - \left(L_i + \delta_1\tau v_i(k) + \delta_2\tau(v_i(k) - v_{i-1}(k))\right)$$

**Proof:** To prove the sequential feasibility of the constraint set $S_1(u_i(k))$ constituted of Equations (2)-(4) and (6)-(8), we need to find a non-empty control input profile $\mathbb{S}_1(u_i(k))$ at step $k$ for $k \in \mathbb{Z}_+$ that makes the constraint set $S_1(u_i(k+1))$ feasible, given that $S_1(u_i(k))$ is feasible. Namely, with feasible state at any step $k$, the MPC can have a feasible control input at step $k$ leading to a feasible state at step $k+1$. Below we provide the technical details.

We first reformulate the speed limit constraint at step $k+1$ according to Equations (3), (4) and (7). And then we find its corresponding feasible control input set $u_i(k) \in \left[\underline{a_{i,v}}, \overline{a_{i,v}}\right]$ as follows in Equation (16).

$$\begin{aligned}
v_{min} &\leq v_i(k+1) \leq v_{max} \\
&\Leftrightarrow v_{min} \leq v_i(k) + \tau\big(u_i(k) - \Delta u_i(k)\big) = v_{max} \\
&\Leftrightarrow v_{min} \leq (1-\tau\varepsilon_i)v_i(k) + \tau(1-\eta_i)u_i(k) + \tau\eta_i u_i(k-1) \leq v_{max} \\
&\Leftrightarrow u_i(k) \in \left[\underline{a_{i,v}}, \overline{a_{i,v}}\right],
\end{aligned} \quad (16)$$

where the lower bound $\underline{a_{i,v}}$ and upper bound $\overline{a_{i,v}}$ are given in Equation (17).

$$\begin{aligned}
\underline{a_{i,v}} &= \frac{v_{min} - (1-\tau\varepsilon_i)v_i(k) - \tau\eta_i u_i(k-1)}{\tau(1-\eta_i)} \\
\overline{a_{i,v}} &= \frac{v_{max} - (1-\tau\varepsilon_i)v_i(k) - \tau\eta_i u_i(k-1)}{\tau(1-\eta_i)}
\end{aligned} \quad (17)$$

Similarly, we reformulate the safe distance constraints at step $k+1$ according to Equations (2)-(4) and (8). For discussion convenience, we use $g_i(k)$ to represent the safe distance constraint in Equation (8) for CAV $i$ at step $k$ and denote $\underline{u_i(k)} = u_i(k) - \Delta u_i(k)$. Then, we present the mathematical derivations below in Equation (18).

$$\begin{aligned}
g_i(k+1) &= x_{i-1}(k+1) - x_i(k+1) \\
&\quad - \Big(L_i + \delta_1 \tau v_i(k+1) + \delta_2 \tau\big(v_i(k+1) - v_{i-1}(k+1)\big)\Big) \\
&= x_{i-1}(k) - x_i(k) + \tau\big(v_{i-1}(k) - v_i(k)\big) + \frac{\tau^2}{2}\Big(\underline{u_{i-1}(k)} - \underline{u_i(k)}\Big) \\
&\quad - \Big(L_i + \delta_1 \tau v_i(k) + \delta_2 \tau\big(v_i(k) - v_{i-1}(k)\big)\Big) - \delta_1 \tau^2 \underline{u_i(k)} \\
&\quad - \delta_2 \tau^2 \Big(\underline{u_i(k)} - \underline{u_{i-1}(k)}\Big) \\
&= g_i(k) + \tau\big(v_{i-1}(k) - v_i(k)\big) + \tau^2\left(\delta_2 + \frac{1}{2}\right)\underline{u_{i-1}(k)} - \tau^2\left(\delta_1 + \delta_2 + \frac{1}{2}\right)\underline{u_i(k)} \\
&= g_i(k) + \tau\big(v_{i-1}(k+1) - v_i(k)\big) + \tau^2\left(\delta_2 - \frac{1}{2}\right)\underline{u_{i-1}(k)} \\
&\quad - \tau^2\left(\delta_1 + \delta_2 + \frac{1}{2}\right)\underline{u_i(k)}
\end{aligned} \quad (18)$$

Note that we assume $S_1(u_i(k))$ is feasible, it makes the safe distance constraints in Equation (8) feasible at time step $k$. Mathematically, $g_i(k) \geq 0$. To make the safe distance constraints keep feasible at next time step $k+1$ (i.e., $g_i(k+1) \geq 0$), we should have the following control input requirement ($u_i(k) \leq \overline{a_{i,d}}$) in Equation (19) based upon Equations (4) and (18).

$$u_i(k) \leq \overline{a_{i,d}} = \varepsilon_i v_i(k) + \eta_i u_i(k) - \eta_i u_i(k-1) \\
+ \frac{g_i(k) + \tau\big(v_{i-1}(k+1) - v_i(k)\big)}{\tau^2\left(\delta_1 + \delta_2 + \frac{1}{2}\right)} + \frac{\left(\delta_2 - \frac{1}{2}\right)\underline{u_{i-1}(k)}}{\delta_1 + \delta_2 + \frac{1}{2}} \quad (19)$$

Recall that we define $\delta_1 \geq 1$, $\delta_2 \geq 0$ for Equation (8). Without loss of generalizability, we pick $\delta_1 \geq \max\left\{\frac{v_{min} - v_{max}}{\tau(a_{min,i} - \varepsilon v_{min})} - 1, 1\right\}$, $\delta_2 = \frac{1}{2}$ to show the sequential feasibility[1]. By plugging in $\delta_2 = \frac{1}{2}$, we can remove the term with $\underline{u_{i-1}(k)}$ and simplify Equation (19) to Equation (20) below.

---

[1] Please note that the selection of the parameters here is to ensure feasibility rigorously. It is not necessarily the best choice for the implementation.

$$\overline{a_{i,d}} = \varepsilon_i v_i(k) + \eta_i u_i(k) - \eta_i u_i(k-1) + \frac{g_i(k) + \tau(v_{i-1}(k+1) - v_i(k))}{\tau^2(\delta_1 + 1)} \tag{20}$$

Wrapping the Equations (6), (16), (19) and (20), we have the following solution set $\mathbb{S}_1(u_i(k))$ in Equation (21) that makes the constraints set $\mathcal{S}_1(u_i(k+1))$ feasible given that $\mathcal{S}_1(u_i(k))$ is feasible.

$$u_i(k) \in \mathbb{S}_1(u_i(k)) = \left[\max\{a_{min,i}, \underline{a_{i,v}}\}, \min\{a_{max,i}, \overline{a_{i,v}}, \overline{a_{i,d}}\}\right] \tag{21}$$

We next show $\mathbb{S}_1(u_i(k))$ in Equation (21) is non-empty. To do that, it suffices to show $\max\{a_{min,i}, \underline{a_{i,v}}\} \leq \min\{a_{max,i}, \overline{a_{i,v}}, \overline{a_{i,d}}\}$. More specifically, we need to prove the following six inequalities hold (i) $a_{max,i} \geq a_{min,i}$, (ii) $a_{max,i} \geq \underline{a_{i,v}}$, (iii) $\overline{a_{i,v}} \geq a_{min,i}$, (iv) $\overline{a_{i,v}} \geq \underline{a_{i,v}}$, (v) $\overline{a_{i,d}} \geq a_{min,i}$, (vi) $\overline{a_{i,d}} \geq \underline{a_{i,v}}$. It is obvious that (i) $a_{max,i} \geq a_{min,i}$ and (iv) $\overline{a_{i,v}} \geq \underline{a_{i,v}}$ hold according to Equations (6) and (17). Below we sequentially show the inequalities (ii), (iii), (v), and (vi) are satisfied in the Equations (22)-(25) when $\delta_1 \geq \max\left\{\frac{v_{min} - v_{max}}{\tau(a_{min,i} - \varepsilon v_{min})} - 1, 1\right\} \geq \frac{v_{min} - v_{max}}{\tau(a_{min,i} - \varepsilon v_{min})} - 1$ and $\delta_2 = \frac{1}{2}$.

Specifically, we confirm inequality (ii) $a_{max,i} \geq \underline{a_{i,v}}$ holds by the derivations given in Equation (22).

$$\begin{aligned} a_{max,i} - \underline{a_{i,v}} &= a_{max,i} - \frac{v_{min} - (1 - \tau\varepsilon_i)v_i(k) - \tau\eta_i u_i(k-1)}{\tau(1-\eta_i)} \\ &\geq a_{max,i} + \frac{-v_{min} + (1-\tau\varepsilon_i)v_{min} + \tau\eta_i a_{min,i}}{\tau(1-\eta_i)} \\ &= a_{max,i} - \frac{\varepsilon_i v_{min} - \eta_i a_{min,i}}{(1-\eta_i)} > 0 \end{aligned} \tag{22}$$

We prove inequality (iii) $\overline{a_{i,v}} \geq a_{min,i}$ by the mathematical process in Equation (23).

$$\begin{aligned} \overline{a_{i,v}} - a_{min,i} &= \frac{v_{max} - (1-\tau\varepsilon_i)v_i(k) - \tau\eta_i u_i(k-1)}{\tau(1-\eta_i)} - a_{min,i} \\ &\geq \frac{v_{max} - (1-\tau\varepsilon_i)v_{max} - \tau\eta_i a_{max,i}}{\tau(1-\eta_i)} - a_{min,i} \\ &= \frac{\varepsilon_i v_{max} - \eta_i a_{max,i}}{(1-\eta_i)} - a_{min,i} > 0 \end{aligned} \tag{23}$$

To ensure inequality (v) $\overline{a_{i,d}} \geq a_{min,i}$, we develop the mathematical process in Equation (24).

$$\overline{a_{i,d}} - a_{min,i} = \frac{\varepsilon_i v_i(k) - \eta_i u_i(k-1)}{(1-\eta_i)} + \frac{g_i(k) + \tau(v_{i-1}(k+1) - v_i(k))}{\tau^2(\delta_1+1)(1-\eta_i)}$$
$$- a_{min,i}$$
$$\geq \frac{1}{(1-\eta_i)}\left[\frac{v_{i-1}(k+1) - v_i(k)}{\tau(\delta_1+1)} - (1-\eta_i)a_{min,i}\right.$$
$$\left. + (\varepsilon_i v_i(k) - \eta_i u_i(k-1))\right] + \frac{g_i(k)}{\tau^2(\delta_1+1)(1-\eta_i)}$$
$$\geq \frac{1}{(1-\eta_i)}\left[\frac{v_{min} - v_{max}}{\tau(\delta_1+1)} - (1-\eta_i)a_{min,i}\right.$$
$$\left. + (\varepsilon_i v_{min} - \eta_i a_{min,i})\right]$$
$$= \frac{1}{(1-\eta_i)}\left[\frac{v_{min} - v_{max}}{\tau(\delta_1+1)} - (a_{min,i} - \varepsilon_i v_{min})\right] \quad (24)$$

By choosing a feasible $\delta_1 \geq \frac{v_{min}-v_{max}}{\tau(a_{min,i}-\varepsilon v_{min})} - 1$, we have $\frac{1}{(1-\eta_i)}\left[\frac{v_{min}-v_{max}}{\tau(\delta_1+1)} - (a_{min,i}-\varepsilon_i v_{min})\right] \geq 0$. Consequently, we confirm inequality (v) $\overline{a_{i,d}} - a_{min,i} \geq 0$ in Equation (24).

Last, we confirm inequality (v) $\overline{a_{i,d}} \geq a_{min,i}$ by the derivation below in Equation (25).

$$\overline{a_{i,d}} - \underline{a_{i,v}} = \left[\frac{\varepsilon_i v_i(k) - \eta_i u_i(k-1)}{(1-\eta_i)} + \frac{g_i(k) + \tau(v_{i-1}(k+1) - v_i(k))}{\tau^2(\delta_1+1)(1-\eta_i)}\right]$$
$$- \frac{v_{min} - (1-\tau\varepsilon_i)v_i(k) - \tau\eta_i u_i(k-1)}{\tau(1-\eta_i)}$$
$$> \frac{1}{\tau(1-\eta_i)}\left[\frac{v_{i-1}(k+1) - v_i(k)}{\delta_1+1} - (v_{min} - v_i(k))\right]$$
$$= \frac{v_{i-1}(k+1) - v_{min} + \delta_1(v_i(k) - v_{min})}{\tau(1-\eta_i)(\delta_1+1)} \geq 0 \quad (25)$$

Wrapping the results above, we prove the sequential feasibility of the constraints $\mathcal{S}_1(u(k))$, with which we conclude **Lemma 1**. ∎

### 3.2. MPC terminal constraints feasibility

This section proves the feasibility of the terminal constraints set $\mathcal{S}_2(u(p), Z(P))$ in Equation (14). Namely, this study wants to find feasible control inputs that satisfies both $\mathcal{S}_1(u(p))$ and $\mathcal{S}_2(u(p), Z(P))$ for $p \in P$. To do that, this study wants to find a lower bound $\underline{P_{\mathbb{E}n}}$ of the prediction horizon $P$ to ensure the feasibility, given that the platoon is under a general scenario $\mathbb{E}_n$. This indicates that when $P \geq \underline{P_{\mathbb{E}n}}$, the feasibility of the constraints sets $\mathcal{S}_1(u(p))$ and $\mathcal{S}_2(u(p), Z(P))$ is ensured so that the feasibility of the problem is proved. We start with a simple case $\mathbb{E}_1$ where only one CAV $i = 1$ in the platoon follows the leading HDV. We investigate the prediction horizon's lower bound $\underline{P_{\mathbb{E}1}}$ to ensure the feasibility of this case $\mathbb{E}_1$ in **Lemma 2**-**Lemma 3**. Without loss of generality, the leading HDV is assumed to drive at a constant speed $v_0$ until the CAV $i = 1$ achieves the steady state $z_1(\rho) = 0, z_1'(\rho) = 0$ under such case $\mathbb{E}_1$. Then we extend the results of the case $\mathbb{E}_1$ in **Lemma 2**-**Lemma 3** to a more general case $\mathbb{E}_n$ with $n$ CAVs in the platoon in **Theorem 1**.

Below this study first illustrates the main idea of the **Lemma 2**-**Lemma 3**. To find the lower bound $\underline{P_{\mathbb{E}1}}$ ensuring the feasibility of the case $\mathbb{E}_1$, we define a special scenario $(E_1)$. **Lemma 2** proves that the case $\mathbb{E}_1$ can be converted into the special scenario $(E_1)$ in limiting time steps. Then **Lemma 3** further proves that there exists a lower bound $P_{E1}$ to ensure the feasibility of the special scenario $(E_1)$. Below we formally define some mathematical notations and the three scenarios $\mathbb{E}_n, \mathbb{E}_1, E_1$.

The HDV's speed is constant and denoted by $v_0$, which satisfies $v_{min} \leq v_0 \leq v_{max}$. The CAV $i = 1$'s speed at time step $p$ is denoted by $v_1(p) = (1 + \delta_p)v_0$, where $\delta_p$ is the discrepancy coefficient between CAV $i = 1$ and the HDV at step $p$. Then according to the vehicle dynamics in Equations (2)-(3) without considering control uncertainties, we have

$$v_0 \delta_{p+1} = v_0 \delta_p + \tau u_1(p) \tag{26}$$

When the leading HDV's speed is constant at $v_0$, the desired speed for the following CAVs are all $v_0$ and the corresponding desired intervehicle spacing is denoted by $s_0 = d(v_0, v_0)$. In particular, the inter-vehicle spacing between CAV $i = 1$ and the HDV is denoted by $s_1(p) = x_1(p) - x_0(p)$. Accordingly, we define the following three scenarios: $\mathbb{E}_n, \mathbb{E}_1, E_1$.

$\mathbb{E}_n$: A platoon with $n$ CAVs follow a leading HDV. All CAVs' initial states satisfy the speed and safe spacing constraints in Equations (6)-(8). Mathematically, for $i \in I$, $v_{min} \leq v_i(0) \leq v_{max}$, $x_i(0) - x_{i-1}(0) \geq d_i(v_i(0), v_{i-1}(0))$.

$\mathbb{E}_1$: A single CAV $i = 1$ follows the leading HDV. $v_{min} \leq v_1(0) \leq v_{max}$, $s_1(0) = x_1(0) - x_0(0) \geq d_1(v_1(0), v_0)$.

$E_1$: A single CAV $i = 1$ follows the leading HDV. The initial speed $v_1(0)$ of the CAV $i = 1$ is larger than the HDV's constant speed $v_0$. Besides, the inter-vehicle spacing between CAV $i = 1$ and the leading HDV is just equal to the safe distance formula in Equation (8). Mathematically, $v_0 \leq v_1(0) \leq v_{max}$, $s_1(0) = x_1(0) - x_0(0) = d_1(v_1(0), v_0)$.

Then under the scenario $(E_1)$, we have $\delta_0 \geq 0$. Based on the three scenarios defined above, below we prove **Lemma 2**.

**Lemma 2.** There exist feasible control inputs $u \in \mathbb{S}_1(u_i(k-1))$ to convert scenario $(\mathbb{E}_1)$ to scenario $(E_1)$ in $\rho_t = \left\lceil \frac{(s_1(0) - d_1(v_{max}, v_0))_+}{\tau v_{max}} \right\rceil + \left\lceil \frac{v_{max} - v_1}{\tau a_{max}} \right\rceil + \left\lceil \frac{(s_0 - s_1(0))_+}{v_0 - v_1(0)} \right\rceil$ time steps.

**Proof:**
This study realized there exists three situations for the scenario $(\mathbb{E}_1)$.
**S(i).** The CAV $i = 1$'s speed is larger than the HDV's speed. However, the inter-vehicle spacing is larger than the safe distance $d_1(v_1, v_0)$. Mathematically, $v_0 \leq v_1(0) \leq v_{max}$, $s_1(0) > d_1(v_1(0), v_0)$.
**S(ii).** The CAV $i = 1$'s speed is no larger than the HDV's speed. And the inter-vehicle spacing is larger than the desired spacing $s_0 = d_1(v_0, v_0)$. Mathematically, $v_{min} \leq v_1(0) < v_0$, $s_1(0) \geq s_0$.
**S(iii).** The CAV $i = 1$'s speed is no larger than the HDV's speed. The inter-vehicle spacing is larger than the safe distance $d_1(v_1, v_0)$, but smaller than the desired spacing $s_0 = d_1(v_0, v_0)$. Mathematically, $v_{min} \leq v_1(0) < v_0$, $d_1(v_1, v_0) \leq s_1(0) \leq s_0$
**The First Situation S(i):**

We first analyze the first situation **S(i)**. Mainly, we develop a strategy $\mathfrak{s}_1$ to make scenario $(\mathbb{E}_1)$ under **S(i)** convert to scenario $(E_1)$ within $t'$ time steps. Mathematically, $v_0 \leq v_1(t') \leq v_{max}, s_1(t') = d_1(v_1(t'), v_0)$. The strategy $\mathfrak{s}_1$ is presented below.

$\mathfrak{s}_1$: Let CAV $i = 1$ accelerate at $a_{max}$ for $\min\left\{\left\lceil \frac{v_{max} - v_1(0)}{\tau a_{max}} \right\rceil, t\right\}$ time steps until the safe distance bound or speed bound will be violated at step $\min\left\{\left\lceil \frac{v_{max} - v_1(0)}{\tau a_{max}} \right\rceil, t\right\} + 1$. Where $t$ is uniquely determined by the following inequalities in Equations (27) and (28).

$$s_1(t) = s_1(0) - \frac{t\tau}{2}(2v_0 + t\tau a_{max}) \geq d_1(v_1 + t\tau a_{max}, v_0) \tag{27}$$

$$s_1(t+1) = s_1(0) - \frac{(t+1)\tau}{2}(2v_0 + (t+1)\tau a_{max}) < d_1(v_1 + (t+1)\tau a_{max}, v_0) \tag{28}$$

Under the strategy $\mathfrak{s}_1$, if $\left|\frac{v_{max}-v_1(0)}{\tau a_{max}}\right| \geq t$, then it will take at most $t+1$ time steps to transit to scenario $(E_1)$. Else if $\left|\frac{v_{max}-v_1(0)}{\tau a_{max}}\right| < t$, then CAV $i=1$ will first accelerate for $\left\lceil\frac{v_{max}-v_1}{\tau a_{max}}\right\rceil$ time steps until its speed reaches $v_{max}$, then it will keep its speed at $v_{max}$ until it reaches the safe distance constraints bound. Accordingly, this process totally takes at most $\left\lceil\frac{s_1(0)-d_1(v_{max},v_0)}{\tau v_{max}}\right\rceil + \left\lceil\frac{v_{max}-v_1}{\tau a_{max}}\right\rceil$ time steps to transit to scenario (E).

In summary, it takes at most $\left\lceil\frac{(s_1(0)-d_1(v_{max},v_0))_+}{\tau v_{max}}\right\rceil + \left\lceil\frac{v_{max}-v_1}{\tau a_{max}}\right\rceil$ time steps to transit from scenario $(\mathbb{E}_1)$ to scenario $(E_1)$ under the first situation **S(i)**.

**The Second Situation S(ii):**

Then we discuss the second situation S(ii), where $s_1(0) \geq s_0 = d_1(v_0, v_0)$. This allows us to follow the same strategy $\mathfrak{s}_1$ in **S(i)** to make $s_1(t') = d_1(v_1(t'), v_0)$. Besides, it is natural to have $v_0 \leq v_1(t') \leq v_{max}$. Below we present a short proof about $v_0 \leq v_1(t') \leq v_{max}$.
Following the strategy $\mathfrak{s}_1$, there exists a time $t'' < t'$ that $v_1(t'') = v_0$ and $s_1(t'') > s_1(0) \geq s_0$. Since there is no deceleration action in the strategy $\mathfrak{s}_1$ and $t' \geq t''$, $v_1(t') > v_1(t'') = v_0$.

In summary, it takes at most $\left\lceil\frac{(s_1(0)-d_1(v_{max},v_0))_+}{\tau v_{max}}\right\rceil + \left\lceil\frac{v_{max}-v_1}{\tau a_{max}}\right\rceil$ time steps to transit from scenario $(\mathbb{E}_1)$ to scenario $(E_1)$ for the second situation **S(ii)**.

**The Third Situation S(iii):**

Finally, we investigate the third situation **S(iii)**, where $d_1(v_1(0), v_0) \leq s_1(0) < s_0$. To do that, we first make CAV $i=1$ keep the speed at $v_1(0)$ for $\left\lceil\frac{s_0-s_1(0)}{v_0-v_1(0)}\right\rceil$ time steps. Then, $v\left(\left\lceil\frac{s_0-s_1(0)}{v_0-v_1(0)}\right\rceil\right) = v_0, s_1\left(\left\lceil\frac{s_0-s_1(0)}{v_0-v_1(0)}\right\rceil\right) \geq s_0$. Consequently, **S(iii)** is converted into **S(ii)** by taking extra $\left\lceil\frac{s_0-s_1(0)}{v_0-v_1(0)}\right\rceil$ time steps.

In summary, it takes at most $\left\lceil\frac{(s_1(0)-d_1(v_{max},v_0))_+}{\tau v_{max}}\right\rceil + \left\lceil\frac{v_{max}-v_1}{\tau a_{max}}\right\rceil + \left\lceil\frac{s_0-s_1(0)}{v_0-v_1(0)}\right\rceil$ time steps to transit from scenario $(\mathbb{E}_1)$ to scenario $(E_1)$ for the second situation **S(iii)**.

Wrapping above, it takes at most $\rho_t = \left\lceil\frac{(s_1(0)-d_1(v_{max},v_0))_+}{\tau v_{max}}\right\rceil + \left\lceil\frac{v_{max}-v_1}{\tau a_{max}}\right\rceil + \left\lceil\frac{(s_0-s_1(0))_+}{v_0-v_1(0)}\right\rceil$ time steps to convert the scenario $(\mathbb{E}_1)$ to scenario $(E_1)$. #####

**Lemma 3.** Under the scenario $(E_1)$, there exists feasible control input set $\mathfrak{s}_1(u_1(p)) = \left[\max\{a_{min}, a_{i,v}\}, \frac{v_0\delta_p}{\tau}(D_p - 1)\right] \in \mathbb{S}_1(u_i(p))$. There exists a feasible control sequence $u_1^*(p) \in \mathfrak{s}_1(u_1(p))$ for $p \in \rho_1$, to make CAV $i=1$ can reach $z_1(\rho_1) = 0, z_1'(\rho_1) > 0$ in $\rho_1 = \left\lceil\log_{D_\infty}\left(1 - \frac{2(1-D_\infty)}{1+D_\infty}\left(\frac{1}{1-D_0} - \sigma\right)\right)\right\rceil$ time steps. Using the feasible control sequence $u_1^{**}(p) = a_{min}$ for $p \in \rho_2$, CAV $i=1$ can reach $z_1'(\rho_2) = 0, z_1(\rho_2) > 0$ in $\rho_2 = \left\lceil\frac{v_0\delta_0}{-a_{min}\tau}\right\rceil$ time steps.
$D_p, D_0, D_\infty$ are defined as
$$\begin{aligned} D_p &= \frac{2v_0 + \delta_p v_0 + 2\tau a_{max} - \tau a_{min}}{2v_0 + \delta_p v_0 + 2\tau a_{max} - 3\tau a_{min}} \\ D_0 &= \frac{2v_0 + \delta_0 v_0 + 2\tau a_{max} - \tau a_{min}}{2v_0 + \delta_0 v_0 + 2\tau a_{max} - 3\tau a_{min}} \end{aligned} \tag{29}$$

$$D_\infty = \frac{2v_0 + 2\tau a_{max} - \tau a_{min}}{2v_0 + 2\tau a_{max} - 3\tau a_{min}}$$

**Proof:** This proof involves three different but related statements. Thus, we will prove them sequentially below in Proof (1), Proof (2) and Proof (3).

**Proof (1).** We first prove the first statement: the proposed control input $s_1(u_1(p))$ is feasible, mathematically $s_1(u_1(p)) = \left[a_{min}, \frac{v_0 \delta_p}{\tau}(D_p - 1)\right] \in \mathbb{S}_1(u_1(p))$. Based upon Equation (15) in **Lemma 1**, it suffices to show that

$$\max\{a_{min}, \underline{a_{1,v}}\} \leq \overline{u_1}(p) := \frac{v_0 \delta_p}{\tau}(D_p - 1) \leq \min\{a_{max}, \overline{a_{1,v}}, \overline{a_{1,d}}\} \tag{30}$$

According to **Lemma 3**, the proposed control input can be presented as follows in Equation (31).

$$\overline{u_1}(p) = \frac{v_0 \delta_p}{\tau}(D_p - 1) = \frac{v_0 \delta_p \cdot 2 a_{min}}{\tau\left(\frac{2v_0}{\tau} + \frac{\delta v_0}{\tau} + 2a_{max} - 3a_{min}\right)} \tag{31}$$

According to Equations (15), we have

$$\overline{a_{1,d}} = \frac{3}{2}a_{min} - a_{max} - \frac{v_i(p)}{\tau} + \frac{\sqrt{\Delta_1(p)}}{2} \tag{32}$$

Where

$$\Delta_1(p) = \frac{4(v_0(1 + \delta_p) + 2\tau a_{max} - \tau a_{min})v_0(1 + \delta_p)}{\tau^2} + \frac{-8 a_{min} v_0}{\tau}$$
$$+ 4\left(a_{max} - \frac{3}{2}a_{min}\right)^2 \tag{33}$$
$$\geq 4\left(a_{max} - \frac{3}{2}a_{min} + \frac{v_0}{\tau} + \frac{D\delta_p v_0}{\tau}\right)^2$$

Then based on Equations (32) and (33), we have

$$\overline{a_{1,d}} \geq \frac{v_0 \delta_p}{\tau}(D_p - 1) = u(p) \tag{34}$$

Utilizing the Equation (31), the following relationship in Equation (35) holds.

$$\frac{\overline{u_1}(p)}{v_0 \delta_p} = \frac{2a_{min}}{\tau\left(\frac{2v_0}{\tau} + \frac{\delta v_0}{\tau} + 2a_{max} - 3a_{min}\right)} \in \left(-\frac{1}{\tau}, 0\right) \tag{35}$$

According to Equations (26) and (35), we can derive

$$v_0 \delta_{p+1} = v_0 \delta_p + \tau u_1(p) \in [0, v_0 \delta_p] \tag{36}$$

Equation (36) indicates that once initially $v_{min} - v_0 \leq v_0 \delta_0 \leq v_{max} - v_0$ holds, then by applying the proposed control input in Equation (31), sequentially we have $v_{min} - v_0 \leq v_p \delta_p \leq v_{max} - v_0$. This indicates that $\underline{a_{1,v}} \leq \overline{u_1}(p) \leq \overline{a_{1,v}}$ holds.

According to the Equation (31), we can show $a_{min} \leq \overline{u_1}(p) \leq a_{max}$ in the following Equations (37) and (38)

$$\overline{u_1}(p) = \frac{v_0 \delta_p \cdot 2 a_{min}}{\tau\left(\frac{2v_0}{\tau} + \frac{\delta_p v_0}{\tau} + 2a_{max} - 3a_{min}\right)} = \frac{a_{min}}{\frac{v_0}{\delta_p} + \frac{1}{2} + \frac{2a_{max} - 3a_{min}}{v_0 \delta_p}} > a_{min} \tag{37}$$

$$\overline{u_1}(p) - a_{max} = \frac{(-\delta_p v_0)(a_{max} - 2a_{min}) - 2v_0 a_{max} - (2a_{max} - 3a_{min})a_{max}}{\tau\left(\frac{2v_0}{\tau} + \frac{\delta_p v_0}{\tau} + 2a_{max} - 3a_{min}\right)} \leq 0 \tag{38}$$

Wrapping above and based on the definition of the $\mathbb{S}_1(u_i(p))$ in Equation (15) in **Lemma 1**, we can conclude that the proposed control input $\mathsf{s}_1(u_1(p))$ is feasible, $\mathsf{s}_1(u_1(p)) = \left[a_{min}, \frac{v_0 \delta_p}{\tau}(D_p - 1)\right] \in \mathbb{S}_1(u_1(p))$. ####

**Proof (2).** We next prove the second statement: when $\rho_1 = \left\lceil \log_{D_\infty}\left(1 - \frac{2(1-D_\infty)}{1+D_\infty}\left(\frac{1}{1-D_0} - \sigma\right)\right)\right\rceil$, CAV $i = 1$ can come to the steady state $z_1(\rho_1) = 0, z_1'(\rho_1) > 0$, using a feasible control input sequence $u_1^*(p) \in \mathsf{s}_1(u_1(p))$ for $p \in \rho_1$. It suffices to show that using the feasible control input sequence $\overline{u_1}(p) = \frac{v_0 \delta_p}{\tau}(D_p - 1)$, $z_1(\rho_1) \leq 0, z_1'(\rho_1) > 0$ can be achieved within $\rho_1$ time steps.

According to the definition of $D_p$ in Equation (29) and the relationship $0 < \delta_{p+1} < \delta_p$ derived from the Equation (36), we have

$$D_{p+1} < D_p \tag{39}$$

Based on the Equations (26) and (31), using the feasible control input sequence $\overline{u_1}(p)$, we have

$$v_1(p) = v_0 + v_0 \delta_p = v_0 + v_0 \delta_{p-1} D_{p-1} = \cdots = v_0 + v_0 \delta_0 \prod_{j=0}^{p-1} D_j \tag{40}$$

$$v_1(p) > v_0 + v_0 \delta_0 D_\infty^p \tag{41}$$

Then, when $\rho_1 \geq \left\lceil \log_{D_\infty}\left(1 - \frac{2(1-D_\infty)}{1+D_\infty}\left(\frac{1}{1-D_0} - \sigma\right)\right)\right\rceil$, we can derive $z_1(\rho_1) < 0$ in the following Equation (42).

$$z_1(\rho_1) = s_1(\rho_1) - s_0 = s_1(0) - s_0 + v_0 \rho_1 \tau - \sum_{p=0}^{\rho_1 - 1} \frac{v_1(p) + v_1(p+1)}{2} \tag{42}$$

$$< \tau \delta_0 v_0 - \frac{\delta_0 v_0 (2\tau a_{max} + 2v_0 + \delta_0 v_0)}{2a_{min}} - \sigma - v_0 \delta_0 \tau \left(\frac{1}{2} + D_\infty + D_\infty^2 + \cdots + D_\infty^{p-1} + \frac{D_\infty^{\rho_1}}{2}\right)$$

$$= \tau \delta_0 v_0 - \frac{\delta_0 v_0 (2\tau a_{max} + 2v_0 + \delta_0 v_0)}{2a_{min}} - \sigma - v_0 \delta_0 \tau \left(\frac{1 - D_\infty^{\rho_1}}{1 - D_\infty} + \frac{D_\infty^{\rho_1} - 1}{2}\right)$$

$$< \frac{1}{1 - D_0} - \sigma - v_0 \delta_0 \tau \left(\frac{1 - D_\infty^{\rho_1}}{1 - D_\infty} + \frac{D_\infty^{\rho_1} - 1}{2}\right) \leq 0$$

Besides, according to the Equation (41), we have

$$z_1'(\rho_1) = v_1(\rho_1) - v_0 > v_0 \delta_0 D_\infty^{\rho_1} > 0 \tag{43}$$

Wrapping above, the second statement is proved.

**Proof (3).** We finally prove the third statement: Using the feasible control sequence $u_1^{**}(p)$ for $p \in \rho_2$, CAV $i = 1$ can reach $z_1'(\rho_2) = 0, z_1(\rho_2) > 0$ in $\rho_2 = \left\lceil \frac{v_0 \delta_0}{-a_{min}\tau}\right\rceil$ time steps. It is moted that $v_1(0) = v_0 + v_0 \delta_0 > v_0 > v_{min}$, then $\max\{a_{min}, \underline{a_{1,v}}\} = a_{min}$ until the CAV

$i = 1$ decelerates to $v_0$. It indicates that the CAV $i = 1$ can decelerate at $a_{min}$ until $z_1' = 0$. Under this situation, $z_1'(\rho_2) \leq 0$ is derived in Equation (44).

$$z_1'(\rho_2) = v_0 \delta_0 + a_{min} \tau \rho_2 \leq 0 \tag{44}$$

$$z_1(\rho_2) = s_1(\rho_2) - s_0 = s_1(0) - s_0 + v_0 \rho_1 \tau - \sum_{p=0}^{\rho_1 - 1} \frac{v_1(p) + v_1(p+1)}{2} \tag{45}$$

$$= \tau \delta_0 v_0 - \frac{\delta_0 v_0 (2\tau a_{max} + 2v_0 + \delta_0 v_0)}{2 a_{min}} - \sigma - \frac{v_0 \delta_0 \rho_2}{2} > 0$$

Wrapping above, the third statement is proved. ######

**Lemma 4.** Under the scenario $(E_1)$, there exists a feasible control sequence $u_1^*(p) \in \mathfrak{s}_1(u_1(p))$ for $p \in \rho$, to make CAV $i = 1$ can reach $z_1(\rho) = 0, z_1'(\rho) = 0$, where $\rho = \rho_1 + \rho_2 = \left\lceil \log_{D_\infty} \left( 1 - \frac{2(1-D_\infty)}{1+D_\infty} \left( \frac{1}{1-D_0} - \sigma \right) \right) \right\rceil + \left\lceil \frac{v_0 \delta_0}{-a_{min}\tau} \right\rceil$.

**Proof:**
According to the second statement of the **Lemma 3**, we can similarly conclude that:
(i) there exists a feasible control input sequence $\overline{u_1^*}(p) \in \mathfrak{s}_1(u_1(p))$ that makes $z_1(\rho) = 0, z_1'(\rho) > 0$, since $\rho > \rho_1$. The corresponding feasible speed profile is denoted by $\overline{v_1}(p)$.
Similarly; (ii) there exists a feasible control input sequence $\underline{u_1^*}(p) = \begin{cases} a_{min} & v_1(p) > v_0 \\ 0 & v_1(p) = v_0 \end{cases} \in \mathfrak{s}_1(u_1(p))$ to make $z_1'(\rho) = 0, z_1(\rho) > 0$. Accordingly, CAV $i = 1$'s speed profile is denoted by $\underline{v_1}(p)$.

It is noted that $\underline{u_1^*}(p)$ takes the lower bound of the feasible control input $a_{min}$ until the speed decelerates to $v_0$. Besides, $\overline{v_1}(p) > v_0$ for $\forall p \in \rho$. Hence, we immediately have

$$\overline{v_1}(p) \geq \underline{v_1}(p) \text{ for } \forall p \in \rho. \tag{46}$$

Then according to (i), (ii) and the continuous property of the feasible speed profile $[\underline{v_1}(p), \overline{v_1}(p)]$, we can conclude that there exists a feasible speed profile $v_1^*(p) \in [\underline{v_1}(p), \overline{v_1}(p)]$ to make $z_1(\rho) = 0, z_1'(\rho) = 0$.

According to the definition of scenario $\mathbb{E}_1$, we define scenario $\mathbb{E}_i$ below.
$\mathbb{E}_i$: A single CAV $i$ follows the CAV $i - 1, i - 2, ..., 1$, HDV, which all drive at the constant speed $v_0$. $v_{min} \leq v_i(0) \leq v_{max}$, $s_i(0) = x_i(0) - x_{i-1}(0) \geq d_1(v_i(0), v_0)$. When $P \geq \underline{P_{\mathbb{E}_i}}$, a platoon under scenario $\mathbb{E}_i$ is feasible. $\underline{P_{\mathbb{E}_i}}$ is used to represent the lower bound of the prediction horizon to ensure the feasibility.

**Theorem 1.** When $P \geq \underline{P_{\mathbb{E}1}} = \rho + \rho_t = \left\lceil \frac{(s_1(0) - d_1(v_{max}, v_0))_+}{\tau v_{max}} \right\rceil + \left\lceil \frac{v_{max} - v_1}{\tau a_{max}} \right\rceil + \left\lceil \frac{(s_0 - s_1(0))_+}{v_0 - v_1(0)} \right\rceil + \left\lceil \log_{D_\infty} \left( 1 - \frac{2(1-D_\infty)}{1+D_\infty} \left( \frac{1}{1-D_0} - \sigma \right) \right) \right\rceil + \left\lceil \frac{v_0 \delta_0}{-a_{min}\tau} \right\rceil$, feasibility of a platoon under the $\mathbb{E}_1$ scenario is ensured. When $P \geq \underline{P_{\mathbb{E}n}} = \sum_{i=1}^{n} \underline{P_{\mathbb{E}i}}$, the feasibility of a platoon under the general $\mathbb{E}$ scenario is ensured.

**Proof:**
According to **Lemma 3**, a platoon under the $\mathbb{E}_1$ scenario takes at most $\rho_t$ time steps to transit to $E_1$ scenario. **Lemma 4** indicates that when $P \geq \rho$, the feasibility for a platoon under the $\mathbb{E}_1$ scenario is ensured. Hence, when $P \geq P_{\mathbb{E}1} = \rho + \rho_t$, the model describing a platoon under the $\mathbb{E}_1$ scenario is feasible.

For a platoon under general $\mathbb{E}$ scenario, there is a naïve strategy $\mathcal{S}$: we can sequentially let CAV $i = 1,2,\ldots,n$ apply feasible control inputs to reach the steady states $z_i(P_{\mathbb{E}i}) = 0, z_i'(P_{\mathbb{E}i}) = 0$. Consequently, when $P \geq P_{\mathbb{E}n} = \sum_{i=1}^{n} P_{\mathbb{E}i}$, the feasibility of a platoon under the general $\mathbb{E}$ scenario is ensured. It is noted that the mathematical representation $P_{\mathbb{E}i}$ shares the same structure with $P_{\mathbb{E}1}$. However, it is likely that CAV $i$'s lower bound $P_{\mathbb{E}i}$ is different from CAV 1's lower bound $P_{\mathbb{E}1}$ when $i \neq 1$. That is because CAV $i$ and CAV 1 have different acceleration/deceleration limits (e.g., $a_{max,i}/a_{min,i} \neq a_{max,1}/a_{min,1}$). If $a_{max,i}/a_{min,i} = a_{max,1}/a_{min,1}$, $P_{\mathbb{E}i} = P_{\mathbb{E}1}$.  ######

**Remark 1.** It is noted that $P_{\mathbb{E}n} = \sum_{i=1}^{n} P_{\mathbb{E}i}$ in **Theorem 1** is a conservative lower bound under the scenario $\mathbb{E}$. It is because the strategy $\mathcal{S}$ makes the CAVs sequentially adjust their speeds and inter-vehicle spacings. Practically, all the CAVs can simultaneously adjust the speeds and spacings to quickly reach the steady states. Hence, the practical lower bound can be taken as $P_{\mathbb{E}n}^* = \max_i\{P_{\mathbb{E}i}\} + \pi$, where $\pi \geq 0$ represents the extra time costs resulted from the CAVs' simultaneous actions. However, it is very hard or even impossible to quantitatively analyze the value of $\pi$.

Wrapping above, we can have a rough estimation about the feasible value of lower bound $P_{\mathbb{E}n}^* \in \left[\max_i\{P_{\mathbb{E}i}\}, \sum_{i=1}^{n} P_{\mathbb{E}i}\right]$. Based on that estimation, we can choose the prediction horizon $P$ as $P = \lambda \max_i\{P_{\mathbb{E}i}\} + (1-\lambda)\sum_{i=1}^{n} P_{\mathbb{E}i}, 0 \leq \lambda \leq 1$. In practice, when $\lambda = 0$, the prediction horizon $P = \sum_{i=1}^{n} P_{\mathbb{E}i}$ leads to too large computation loads, we may carefully increase the $\lambda$ to decrease the computation loads. Besides, in real applications, the terminal constraint in Equation (14) may be neglected.

4. **Conclusion**

This study investigates the sequential and terminal constraints feasibility of vehicle platooning MPC model. It is a supplement to our existing paper about to be published.